\documentclass[12pt]{article}
\usepackage{amsmath}\usepackage{amssymb}\usepackage{amscd}

\usepackage[active]{srcltx}
\newcounter{theorem}\makeatletter
\newcounter{parnum}

\newtheorem{theorem}{Theorem}\newtheorem{prop}[theorem]{Proposition}
\newtheorem{lem}[theorem]{Lemma}\newtheorem{def-lem}[theorem]{Definition-Lemma}
\newtheorem{def-prop}[theorem]{Definition-Proposition}
\newtheorem{cor}[theorem]{Corollary}

\def\Ar{{\rm A}}\def\Ab{\bar{\rm A}}
\def\ad{\operatorname{ad}}

\def\ba{\begin{eqnarray}}\def\ea{\end{eqnarray}}\def\be{\begin{equation}}\def\ee{\end{equation}}

\def\f{\mathfrak{k}}\def\h{\mathfrak{h}}\def\g{\mathfrak{g}}\def\sl{\mathfrak{sl}}\def\gl{\mathfrak{gl}}
\def\Ib{\tts\overline{\nns\rm I}_+\tts}\def\Jb{\,\overline{\!\rm I}_-\ts}\def\k{\mathfrak{k}}
\def\mult{\diamond}\def\nminus{\mathfrak{n}_-}
\def\nplus{\mathfrak{n}_+}
\def\nns{\hskip-.5pt}\def\p{\mathfrak{p}}
\def\pp{\widetilde{\phantom{\,}p\phantom{\,}\phantom{\ }}\phantom{\!}\phantom{\!}\phantom{\!}}
\let\rf=\r

\def\ts{\hskip1pt}\def\tts{\hskip.5pt}
\def\U{\mathrm{U}}\def\Uh{\overline{\U}(\h)}
\def\Z{\mathrm{Z}}
\def\Zgk{\Z(\g,\k)}
\def\Zgkk{{\Z}(\g',\k)}

\newcommand{\ZZ}{{\mathbb Z}}

\renewcommand{\theequation}{{\thesection}.{\arabic{equation}}}
\newcommand\qedsymbol{\hbox{\rlap{$\sqcap$}$\sqcup$}}\newcommand\qed{\relax\ifmmode\else\unskip\quad\fi\qedsymbol}

\def\gr{\operatorname{gr}}
\newcommand{\filta}[1]{\Z(\g,\k)^{(#1)}}
\newcommand{\filtc}[1]{\Z(\g,\k)^{(\mult #1)}}
\newcommand{\filtb}[1]{\Z'(\g,\k)_{(#1)}}
\newcommand{\filtd}[1]{\Uh_{(#1)}}
\def\degreeb{\partial} 
\begin{document}

\begin{center}
{\Large\bf Zero divisors in reduction algebras}

\bigskip
\medskip
{\bf  S. Khoroshkin$^{\circ\diamond}$ \ \ and \, O. Ogievetsky$^{\star\sharp}$\footnote{On leave of absence
from P.N. Lebedev Physical Institute, Theoretical Department, Leninsky prospekt 53, 119991 Moscow, Russia}
}\medskip\\
$^\circ${\it Institute of Theoretical and Experimental Physics, 117218
Moscow, Russia}\\
$\diamond$ {\it Higher School of Economics, Myasnitskaya 20, 101000, Moscow, Russia} \smallskip\\
$^\star${\it Centre de Physique Th\'eorique\footnote{Unit\'e Mixte de Recherche
(UMR 6207) du CNRS et des Universit\'es Aix--Marseille I,
Aix--Marseille II et du Sud Toulon -- Var; laboratoire affili\'e \`a la FRUMAM (FR 2291)},
Luminy, 13288 Marseille, France}\\
$\sharp$ {\it French-Russian Poncelet laboratory, UMI 2615 du CNRS\\
Independent University of Moscow, 11 B. Vlasievski per., 119002 Moscow, Russia
}
\end{center}

\smallskip
\setcounter{equation}{0}
\begin{abstract}\noindent
We establish the absence of zero divisors in the reduction algebra of a Lie algebra
$\g$ with respect to its reductive Lie sub-algebra $\f$. The class of reduction algebras include the
Lie algebras (they arise when $\f$ is trivial) and the Gelfand--Kirillov conjecture extends naturally
to the reduction algebras. We formulate the conjecture for the diagonal reduction algebras
of $\sl$ type and verify it on a simplest example.
\end{abstract}

 \setcounter{equation}{0}
\section{Preliminaries}\label{secredalg}
Let $\k$ be a reductive Lie subalgebra of a Lie algebra $\g$; that
is, the adjoint action of $\k$ on $\g$ is completely reducible (in particular,  
$\k$ is reductive). Fix a triangular
decomposition of the Lie algebra $\k$,
\be\label{intro1} \k=\nminus+\h+\nplus\ .\ee Denote by
$\Delta_+$ and $\Delta_-$ the sets of positive and negative roots in
the root system $\Delta=\Delta_+\cup \Delta_-$ of $\k$. {}For each
root $\alpha\in\Delta$ let $h_\alpha=\alpha^\vee\in\h$ be the
corresponding coroot vector. Denote by $\Uh$ the ring of fractions
of the commutative algebra $\U(\h)$ relative to the set of
denominators
\begin{equation}
\label{M2} \{\,h_\alpha+l\ |\ \alpha\in\Delta\ts,\ l\in\ZZ\,\ts\}\,.
\end{equation}
The elements of this ring can also be regarded as rational functions
on the vector space $\h^\ast\ts$. The elements of
$\U(\h)\subset\,{\!\Uh\!\!\!}\,\,\,$ are then regarded
as polynomial functions on $\h^\ast\ts$. Let
$\overline{\U}(\k)\,\,\,\subset \Ab= \overline{\U}(\g)\,\,\,$ be the rings of
fractions of the algebras $\U(\k)$ and $\Ar=\U(\g)$ relative to the set of denominators
\eqref{M2}. These rings are well defined, because both $\U(\f)$ and $\U(\g)$ satisfy the Ore
condition relative to \eqref{M2}; we give a short proof in the second part of Appendix. 

Define $\Z(\g,\k)$ to be the double coset space of $\Ab$ by its left ideal $\Ib:=\Ab\nplus$,
 generated by elements of $\nplus$, and the right ideal $\Jb:=\nminus\Ab$,
generated by elements of $\nminus$, $\Z(\g,\k):=\Ab/(\Ib+\Jb)$.
The space $\Z(\g,\k)$ is an associative algebra with respect to the multiplication map
\begin{equation}\label{not5a}a\mult b:=aP b\ .\end{equation}
Here $P$ is the extremal projector \cite{AST} of the Lie algebra $\k$ corresponding to the triangular
decomposition \rf{intro1}. We call $\Z(\g,\k)$
the {\textit{reduction}} algebra associated to the pair $(\g,\k)$. See \cite{Zh,KO} for details.

Let $\p$ be an $\ad_\k$-invariant complement of $\k$ in $\g$.
Choose a linear basis $\{p_K\}$ of $\p$, $K$ runs through a certain set $\mathfrak{I}$ of indices.
 We assume that
\begin{itemize}
 \item[1) ] each basis vector is a weight vector, that is
\begin{equation}[h,p_K]=\mu_K(h) p_K;
 \label{weir}
\end{equation}

\item[2) ] the set $\mathfrak{I}$ of indices of basis vectors is equipped with a total order $\preceq$, compatible with
the natural order of their weights, that is, if $\mu_K-\mu_L$ is a sum of simple roots of $\k$ with integer
nonnegative coefficients, then necessarily $K\preceq L$.
 \end{itemize}

The reduction algebra  has the following general properties, see \cite{Zh,KO,KO1}:

\begin{itemize}
\item[(i)] $\Zgk$ is
free as a left $\Uh$-module and as a right $\Uh$-module. As a generating
 (over $\Uh$) subspace one can take a projection of the space $\mathrm{S}(\p)$ of
symmetric tensors on $\p$ to $\Zgk$, that is a subspace of $\Zgk$, formed by linear combinations
of images of the powers $p^\nu$, where $p\in\p$ and $\nu\geq 0$. 
Assignments  $\deg (\widetilde{\hspace{.05cm}X\hspace{.05cm}})=l$ for the image
 of any product of $l$ elements from $\p$, $X=p_{K_1}p_{K_2}\cdots p_{K_l}$,
and  $\deg (Y)=0$ for any $Y\in\Uh$ define the structure of a
filtered algebra on $\Zgk$. The subspace $\filta{k}$ of elements of
degree not greater than $k$ is a free left $\Uh$-module and a free
right $\Uh$-module, with a generating subspace formed by linear
combinations of images of the powers $p^\nu$, where $p\in\p$ and
$k\geq\nu\geq 0$. Moreover,
the images of monomials ($\bar{L}$ is understood as the multi-index)
\ba&& p_{\bar{L}}:=p_{L_1}^{n_1}p_{L_2}^{n_2}\cdots p_{L_m}^{n_m}\ ,\quad L_1\prec L_2
\prec\ldots\prec L_m\
,\quad k=n_1+\dots+n_m\ ,\label{inimp}\ea
in $\filta{k}$ are linearly independent over $\Uh$ and their projections to  the quotient
$\filta{k}/ \filta{k-1}$ form a basis of the left $\Uh$-module
$\filta{k}/ \filta{k-1}$. %

\item[(ii)] The algebra $\Zgk$ is the unital associative algebra, generated by $\Uh$ and all
$\pp_L$,
with the weight relations \rf{weir} and the ordering relations
\be\label{not3} \pp_I\mult \pp_J=\displaystyle{\sum_{K,L:K\preceq L}}{\mathrm{B}}_{IJKL}\,
\pp_K\mult\pp_L+\sum_L{\mathrm{C}}_{IJL}\pp_L+ {\mathrm{D}}_{IJ}\ ,\quad I\succ J\ , \ee
where ${\mathrm{B}}_{IJKL}$, ${\mathrm{C}}_{IJL}$  and ${\mathrm{D}}_{IJ}$ are certain
elements of $\Uh$.

\item[(iii)] Assume that $\g$ is finite-dimensional. Then the monomials
\ba&& \pp_{I_1}\mult\pp_{I_2}\mult\cdots\mult\pp_{I_a},\qquad I_1\preceq  I_2\preceq \ldots\preceq I_a\ ,\label{inimpb}\ea
form a basis of the left $\Uh$-module $\Z(\g,\k)$.

\vskip .2cm
Moreover, any expression in $\Z(\g,\k)$ can be written in the ordered form by a repeated
application of \rf{not3} as instructions "replace the left hand side by
the right hand side".
 \end{itemize}

The main goal of this note is to extend this list by the property
\begin{itemize}
 \item[(iv)] The algebra $\Z(\g,\k)$ has no zero divisors.
\end{itemize}
We  discuss the possible generalizations of the  Gelfand-Kirillov conjecture
for reduction algebras; we also  present a proof
 of the property (iii), different from the proof in \cite{KO1}.

\setcounter{equation}{0}
\section{Absence of zero divisors}\label{div0}  Let $\mathfrak{R}^\prec$ be the set of ordering relations \rf{not3}.
The right hand sides of relations in $\mathfrak{R}^\prec$ contain quadratic, linear and degree zero terms;
the algebra $\Z(\g,\k)$ is thus filtered by the degree in the  generators $\pp_L$; 
in the previous section, the members of this filtrations were denoted by $\filta{k}$.

By the standard argument, it is sufficient to prove the
absence of zero divisors in the associated graded algebra.
Note that this graded algebra  is isomorphic to the reduction algebra, related to the pair 
$({\g}',\k)$, where
$\g'$ is the semidirect sum of $\k$ and its module $\p$. In other words, $\g'$ is isomorphic to
$\k\oplus\p$ as a vector space, the Lie brackets between the elements of $\k$ and between elements of $\k$ and of $\p$
 are the same as in $\g$, the other brackets vanish.  
We thus denote the associated graded algebra by  $\Zgkk$;  its generators 
 we shall denote by the same symbols $\pp_L$.

We have also the structure of a filtered space on $\Uh$. The filtration is given by the degree $\degreeb$
of a rational function (the degree is defined to be the difference of the total degrees of
the numerator and denominator as polynomials in several variables).
The subspace $\filtd{k}\subset\Uh$, $k\in \ZZ$, consists of rational functions on $\h^*$ of degree
not greater than $k$. We have
$$ \filtd{k}\subset \filtd{k+1}, \qquad \cap_k \filtd{k}=0,\qquad \cup_k\filtd{k}=\Uh.$$
 Moreover, $h h'\in\filtd{k+l}$, if $h\in\filtd{k}$ è $h'\in\filtd{l}$, that is, $\Uh$
is a filtered ring, and the associated graded quotient $\gr \Uh$ is isomorphic to the ring
  $\breve{\mathrm U}(\h)$ of rational functions on
 $\h^*$ with poles on hyperplanes $h_\gamma=0$,  $\gamma\in\Delta.$

Recall that the algebra $\Zgkk$ is a free left (and right) $\Uh$-module with a basis $\pp_{\bar{L}}$, see \rf{inimp}.
 For any$z=h_1\pp_{\bar{K}_1}+\cdots +h_m\pp_{\bar{K}_m}\in\Zgkk$ with $h_j\in\Uh$ we set 
$\degreeb( z)\leq k$ if $\degreeb( h_j)\leq k$ 
  for all $j$.  This definition does not depend on a choice of a linear basis in $\mathrm{S}(\p)$: 
one can choose instead of $\pp_{\bar{K}}$ the image of an arbitrary basis of symmetric algebra 
of $\p$. Let $\filtb{k}$, where $k\in\ZZ$, be the subspace in $\Zgkk$ formed by elements of degree not greater than 
$k$. We have
$$ \filtb{k}\subset \filtb{k+1}, \qquad \cap_k \filtb{k}=0,\qquad \cup_k\filtb{k}=\Zgkk.$$
\begin{lem}$\hspace{-.2cm}.$
\begin{itemize}
 \item[(i)]  $\Zgkk$ is a filtered algebra with respect to the filtration $\{\filtb{k}\}$. The filtrations on $\Zgkk$ and
 on $\Uh$ are compatible, that is $\Zgkk$ is the filtered module over the filtered ring $\Uh$.
\item[(ii)] The associated graded quotient algebra $\gr \Zgkk$ is isomorphic to the tensor product 
$\breve{\mathrm U}(\h)\otimes \mathrm{S}(\p)$.
\end{itemize}
\end{lem}
{\it Proof}. The nontrivial part of both statements concerns the multiplication structure and immediately follows from the structure
 of the extremal projector $P$: it has the form $$P=1+\sum h_i x_i y_i\ ,$$
where $h_i\in \Uh$, $\partial (h_i)<0$ and $x_i\in \U (\nminus )$, $y_i\in \U (\nplus )$.

Moving in $a\mult b\equiv aPb\mod (\Ib +\Jb )$ the elements $x_i$ to the left through $a$ and $y_i$ to the right 
through $b$ we find that $a\mult b=ab\ +$ terms of lower degree $\partial$. \hfill $\Box$

\vskip .2cm
\begin{cor}$\hspace{-.2cm}.$
 The algebra $\Zgk$ has no zero divizors.
\end{cor}
{\it Proof}. The commutative algebra $\breve{\mathrm U}(\h)\otimes \mathrm{S}(\p)$ clearly has no zero divisors. The absence
of zero divisors in the filtered algebra $\Zgkk$ follows, since its associated graded quotient algebra is 
$\breve{\mathrm U}(\h)\otimes \mathrm{S}(\p)$. Then the absence
of zero divisors in the filtered algebra for $\Zgk$ follows, since its associated graded quotient algebra is $\Zgkk$.
\hfill $\Box$

\setcounter{equation}{0}
\section{Quotient rings of reduction algebras}
Suppose that $\g$ is finite-dimensional.
 The arguments of the previous section imply the noetherian property of reduction algebras. Indeed,
the commutative ring $\breve{\mathrm U}(\h)$ is noetherian as a localization of the noetherian polynomial ring $\U(\h)$. The tensor product 
$\breve{\mathrm U}(\h)\otimes \mathrm{S}(\p)$ is noetherian as well by the Hilbert basis theorem: 
``if $R$ is a left noetherian ring, then the polynomial ring $R[x]$\, is also a left noetherian ring''. The filtered ring is
noetherian if the associated graded quotient is noetherian, see e.g., \cite{GK}. This implies that both $\Zgkk$ and $\Zgk$ are noetherian rings.
 We summarize all statements as
\begin{prop}$\hspace{-.2cm}.$
 Suppose $\g$ is finite-dimensional. Then the reduction algebra $\Zgk$ is a (left) noetherian domain. 
\end{prop}

A (left) noetherian domain is an Ore domain, see \cite{GK}, Lemma 2. We can thus form fields of fractions of reduction 
algebras. It is natural to conjecture that there is a relation between the validity of the Gelfand--Kirillov conjecture
for the field of fractions of $\U (\g )$ and the field of fractions of $\Zgk$ (maybe up to finite extensions
of the centers of the fields of fractions, as in \cite{GK2}). In particular, we conjecture that the Gelfand--Kirillov 
conjecture holds for the diagonal reduction algebra DR$\bigl(sl(n)\bigr)$ (see \cite{KO2} for definitions and notation).  
 More precisely, we conjecture that the field of fractions of DR$\bigl(sl(n)\bigr)$ is isomorphic to the field of fractions
of the algebra $\mathcal{A}_{\frac{n(n-1)}{2}, 2(n-1)}$, where $\mathcal{A}_{k,l}$ is the unital associative 
algebra, generated by $2k+l$ variables $u_1,...,u_k$, $v_1,...,v_k$, $y_1,...,y_l$, with  relations 
$[u_i,v_j]=\delta_{i,j}$, $[u_i,u_j]=[v_i,v_j]=0$,  $[y_\alpha, u_i]=[y_\alpha,v_i]=0$ and $[y_\alpha,y_\beta]=0$. 
Note that the field of fractions
 of $\U(sl(n))$ is isomorphic to  the field of fractions of the algebra $\mathcal{A}_{\frac{n(n-1)}{2}, n-1}$.

\paragraph{Example.} The diagonal reduction algebra DR$\bigl(sl(2)\bigr)$ for the Lie algebra $sl(2)$, see \cite{KO2},
 has generators $z_+,z_-,t$ and $h$ with the defining relations
\be z_+t=tz_+\frac{h+4}{h+2}\ ,\quad z_+z_-=z_-z_+\frac{h(h+3)}{(h+1)(h+2)}-t^2\frac{1}{h}+h\ ,
\quad tz_-=z_-t\frac{h+2}{h}\label{defre1}\ee
and the $h$-weight relations
\be [h,z_+]=2z_+\ ,\qquad [h,t]=0\ ,\qquad [h,z_-]=2z_-\ .\label{defre2}\ee
The central elements are
$$C^{(1)}=(h+2)t\ ,$$
$$\ C^{(2)}=z_-z_+\frac{h+3}{h+2}+t^2\frac{1}{4}+\frac{h(h+4)}{4}\ .$$
The following formulas define  a homomorphism from the algebra DR$\bigl(sl(2)\bigr)$ to a certain 
localization of the algebra ${\mathcal A}_{1,2}$ (with the Weyl variables $x$ and $\frac{d}{dx}$ and commuting variables
$\nu$ and $\zeta$):  
\be z_-\mapsto x^{-1}\ ,\qquad t\mapsto \frac{\nu}{2(E+1)}\ ,\qquad  z_+\mapsto xf(E)\ ,\qquad h\mapsto 2E\ .
\label{iso1}\ee
Here $E$ is the Euler operator,
$$E=x\frac{d}{dx}\ ,$$
and
$$f(E):=-\frac{2(E+1)}{2E+3}\Bigl( E(E+2) +\frac{\nu^2}{16(E+1)^2}+\zeta\Bigr)\ .$$

Now, the following formulas define  a homomorphism from the algebra  ${\mathcal A}_{1,2}$  to a certain localization of DR$\bigl(sl(2)\bigr)$:
\be x\mapsto z_-^{-1}\ ,\qquad \frac{d}{dx}\mapsto \frac{1}{2}z_-h\ ,\qquad \nu\mapsto 2C^{(1)}\ ,\qquad \zeta\mapsto C^{(2)}\ .
\label{iso2}\ee
The extensions to the fields of fractions of the maps %
\rf{iso1} and \rf{iso2} are inverse to each other and establish an isomorphism of the fields of fractions of the algebras  
DR$\bigl(sl(2)\bigr)$ and $\mathcal{A}_{1,2}$.

\bigskip
Moreover, through the homomorphism \rf{iso1} 
the operators of the generators of the algebra DR$\bigl(sl(2)\bigr)$ act naturally on the space 
$V_M$ with the basis $\{ x^{j+M}\,\vert\, j\in\mathbb{Z}\}$; here $M$
is either a parameter from $\mathbb{C}/\mathbb{Z}$ or can be considered as a variable. This construction provides thus a 
two-parameter family of representations of the algebra DR$\bigl(sl(2)\bigr)$.
The family is unique in the following sense.

\vskip .2cm
\begin{prop}$\hspace{-.2cm}.$
 Let $V$ be the space with the basis $\{ v_j\,\vert\, j\in\mathbb{Z}\}$. Assume that the operators
$z_+,z_-,t$ and $h$ act on $V$ by the formulas
$$z_-\colon v_j\mapsto v_{j-1}\ ,\qquad t\colon v_j\mapsto \beta_jv_j\ ,\qquad z_+\colon v_j\mapsto \gamma_jv_{j+1}\ ,\qquad h\colon v_j\mapsto \alpha_jv_j\ $$
with all coefficients $\alpha_j,\beta_j,\gamma_j$ non-vanishing. Then this module is isomorphic to $V_M$ with a non-integer $M$.
\end{prop}

\vskip .2cm\noindent
{\it Proof.} The defining relations \rf{defre1}--\rf{defre2} imply (taking into account that the coefficients are invertible):
$$(\alpha_j+2)\beta_j=\alpha_j\beta_{j-1}\ ,  \quad\alpha_j=\alpha_{j-1}+2\ ,
\quad\gamma_{j-1}=\alpha_j-\frac{\beta_j^2}{\alpha_j}+\frac{\alpha_j(\alpha_j+3)}{(\alpha_j+1)(\alpha_j+2)}\gamma_j\ .$$
The first  two of these relations imply
$$\alpha_j=2j+M\ \ \text{with some $M$ and}\ \ \beta_j=\frac{\nu}{\alpha_j+2}\ \ \text{with some $\nu$}\ .$$
Let
$$\tilde{\gamma}_j:=\dfrac{\alpha_j+3}{\alpha_j+2}\gamma_j\ .$$
Then the recurrence for $\gamma_j$ becomes
$$\tilde{\gamma}_{j-1}=\tilde{\gamma}_j+\frac{\alpha_j+1}{\alpha_j}\Bigl(\alpha_j-\frac{\beta_j^2}{\alpha_j}\Bigr)\ .$$
Substituting the expression for $\beta_j$ and using the identity
$$\frac{4(y+1)}{y^2(y+2)^2}=\frac{1}{y^2}-\frac{1}{(y+2)^2}$$
one easily solves the recurrence for $\gamma_j$ and obtains the assertion.\hfill$\Box$

\renewcommand{\theequation}{{A}.{\arabic{equation}}}
\setcounter{equation}{0}
\section*{Appendix}

\paragraph{1. Proof of the statement (iii), section \ref{secredalg}.} The present proof uses a result about cubic monomials 
from \cite{KO1} and then refers to the diamond, or composition, lemma, see \cite{Bo,Be}. We shall prove the statement (iii) in 
several steps.

Denote by  ${\cal{I}}(\pp_{I_1}\mult\pp_{I_2})$ the right hand side of \rf{not3}. We understand \rf{not3} as
the set of instructions $\pp_{I_1}\mult\pp_{I_2}\leadsto {\cal{I}}(\pp_{I_1}\mult\pp_{I_2})$ ($\leadsto$ stands for "replace") in the free algebra with the weight generators $\pp_I$.

\begin{lem}\hspace{-.2cm}.  Any polynomial, cubic in the generators, acquires an ordered form after a repeated application of instructions \rf{not3}.
\label{lemcub}\end{lem}
The proof is given in \cite{KO1}.

\vskip.2cm
{\it{Proof}} of statement (iii), section \ref{secredalg}.  

\vskip .2cm
(a) Recall the filtration  $\filta{k}$ defined in (i), section \ref{secredalg}.
By the statement (ii), section \ref{secredalg}, the algebra $\Zgk$ 
is generated by $\pp_I$
and, due to the form \rf{not3} of relations, has a filtration by the $\mult$-degree. Let $\filtc{k}$ 
 be the subspace of elements of degree not greater
than $k$ with respect to the product $\mult$. Since $\pp_{I_1}\mult\pp_{I_2}\mult\dots\mult\pp_{I_k}=\pp_{I_1}P\pp_{I_2}P\dots P\pp_{I_k}$, it follows that
$\filtc{k}\subset \filta{k}$.
The opposite inclusion holds as well because the algebra $\Zgk$  
is generated by $\pp_I$. We conclude that the two
filtrations coincide.

Therefore, every
$p_{\hspace{-.03cm}I}\hspace{-.25cm}\widetilde{\hspace{.22cm}p_{\hspace{-.05cm}J}\hspace{.2cm}}\hspace{-.2cm}p_{\hspace{-.03cm}K}$,
$I\preceq J\preceq K$, is in $\filtc{3}$ %
and, by lemma \ref{lemcub}, can be ordered. The cardinalities of the sets
$\{ p_{\hspace{-.03cm}I}\hspace{-.25cm}\widetilde{\hspace{.22cm}p_{\hspace{-.05cm}J}\hspace{.2cm}}\hspace{-.2cm}p_{\hspace{-.03cm}K}
\mid I\preceq J\preceq K\}$ and $\{\pp_I\mult\pp_J\mult\pp_K\mid I\preceq J\preceq K\}$ are equal,
so due to \rf{inimp} the image of the 
set $\{\ \pp_I\mult\pp_J\mult\pp_K\mid I\preceq J\preceq K\}$ is a basis of $\filtc{3}/\filtc{2}$.  
To generalize
this statement to higher degrees, we use the diamond lemma.

\vskip .2cm
(b)  We shortly remind a slightly simplified version of the diamond lemma assertion, following \cite{O}.
Let ${\mathsf{A}}$ be the free associative algebra on letters $\{ x^1,\dots ,x^N\} $. Fix an order on the set $\{ x^1,\dots ,x^N\} $. Write an
element $f\in {\mathsf{A}}$
in the reduced form, as a linear combination of different words. Denote by $\hat{f}$ the {\it highest symbol} of $f$, that is, the lexicographically highest word in
the reduced form of the element $f$.

Let ${\mathsf{B}}$ be a quotient algebra of ${\mathsf{A}}$ by a set ${\mathsf{S}} =\{ r_1,\dots ,r_M\} $ of relations. Every relation $r$ we
write in the form $\hat{r} =\lfloor r\rfloor$, all terms in $\lfloor r\rfloor$ are smaller than $\hat{r}$; we understand it as an instruction to replace $\hat{r}$ by the right
hand side. Taking, if necessary, linear combinations of relations, we always assume that all $\hat{r}$'s are different. Let
$\hat{{\mathsf{S}}}=\{ \hat{r}_1,\dots ,\hat{r}_M\} $.

Given an expression, it might happen that there are different ways of applying instructions to it. This is called {\it ambiguity}. An ambiguity happens iff there are
(maybe after several applications of the instructions) two or more places in which different subwords of the form $\hat{r}$,  $\hat{r}\in \hat{{\mathsf{S}}}$, enter the
expression.

Suppose that all possible sequences of applications of the instructions to every word terminate and lead to one and the same result; so the result
depends on the initial word but not on the way of using the instructions. By construction, the resulting final expression does not contain any subword from
$\hat{{\mathsf{S}}}$. In this situation one says that all ambiguities are {\it resolvable}.

Two types of ambiguities are called {\it minimal}; these are {\it overlaps}
and {\it inclusions}. The overlap ambiguity is of the sort $\hat{r}_i={\mathsf{ab}}$ and $\hat{r}_j={\mathsf{bc}}$ for some words
${\mathsf{a}},{\mathsf{b}},{\mathsf{c}}$ and some $i,j=1,\dots ,M$. The inclusion ambiguity
is of the sort $\hat{r}_i={\mathsf{abc}}$ and $\hat{r}_j={\mathsf{b}}$ for some words
${\mathsf{a}},{\mathsf{b}},{\mathsf{c}}$ and some $i,j=1,\dots ,M$.

The diamond lemma states the equivalence of three assertions:
\begin{itemize}
\item[(1)] all ambiguities are resolvable;
\item[(2)] all minimal ambiguities are resolvable;
\item[(3)] ${\mathsf{B}}$, as a vector space, possesses a basis consisting of images of {\it normal} words - words, which do not have subwords belonging to
$\hat{{\mathsf{S}}}$.
\end{itemize}

We return to our situation (the reduction algebra). The only minimal ambiguities for the instructions \rf{not3} are overlaps
$\pp_{I_1}\mult\pp_{I_2}\mult\pp_{I_3}$ with $I_1\succ I_2\succ I_3$. We can start by either ${\cal{I}}_{12}$ or ${\cal{I}}_{23}$; each time we will arrive, by
lemma \ref{lemcub}, to an ordered expression. By (a) of the present proof, the ordered expressions are linearly independent. Therefore, the two expressions
coincide and the ambiguities are resolvable. The set of normal words is given by \rf{inimpb}, which finishes the proof of the PBW property for the reduction
algebra. \hfill $\qed$

\paragraph{2. Ore conditions.} Let $A$ be an associative algebra and $S$ a multiplicative (that is, multiplicatively closed)
subset in $A$. The Ore conditions, ensuring that one can localize with respect to the "set of denominators" $S$, are
\begin{equation}\label{loou1}\forall\ a\in A,s\in S\ \exists\ \tilde{a}\in A,\tilde{s}\in S\colon\ a\tilde{s}=s\tilde{a}\ ,\end{equation}
or $aS\cap sA\neq\varnothing$, and
\begin{equation}\label{loou2}\text{if}\ sa=0\ \text{for some}\ s\in S,a\in A\ \text{then}\ \exists\ \tilde{s}\in S\colon\ a\tilde{s}=0\ .\end{equation}

Let $S$ be multiplicativly generated by elements $h_{\mu}+k$, $k\in \mathbb{Z}$; we assume that the
algebra $A$ is a sum of its $\mathfrak{h}$-weight components ($\mathfrak{h}$ is the Cartan subalgebra of
$\k$), that is, each element in $A$ can be written as a sum
of  $\mathfrak{h}$-homogeneous elements; recall that an element $x$ is homogeneous of weight $\alpha\in \mathfrak{h}^*$
if $\hat{h}(x)=\alpha(h)x$ for any $h\in \mathfrak{h}$. For an arbitrary
element $h$ we denote by $\hat{h}$  the commutator with $h$, $\hat{h}(x):=[h,x]$.

 \vskip .2cm
We shall verify the Ore conditions in this situation.

\paragraph{Condition (\ref{loou1}).} It is clearly sufficient to check this condition only for multiplicative generators of
the set $S$. Let  $s=h_{\mu}+k$ and let $a$ be an arbitrary element of $A$. Decompose $a$ into a sum of
$\mathfrak{h}$-homogeneous components, $a=a_1+...+a_N$. Then $\hat{h}_{\mu}(a_j)=\mu_ja_j$
with some numbers $\mu_j$ or $a_j (h_{\mu}+k+\mu_j)=(h_{\mu}+k)a_j$ for all $j=1,...,N$. Therefore, for
$\tilde{s}:=(h_{\mu}+k+\mu_1)(h_{\mu}+k+\mu_2)...(h_{\mu}+k+\mu_N)$ the product
$$a\tilde{s}=(a_1+...+a_N) (h_{\mu}+k+\mu_1)(h_{\mu}+k+\mu_2)...(h_{\mu}+k+\mu_N)$$
is divisible from the left by $h_{\mu}+k$, that is, representable in the form $s\tilde{a}$.\hfill$\Box$

\paragraph{Condition (\ref{loou2}).} Again, it is enough to check the condition for generators of $S$ only: if the condition
(\ref{loou2}) holds for

 \vskip .1cm
$\bullet$ $s_1\in S$ and all $\{a\in A\colon s_1a=0\}$
 \vskip .1cm

\noindent and for

 \vskip .1cm
$\bullet$ $s_2\in S$ and all $\{ a\in A\colon s_2a=0\}$
 \vskip .1cm

\noindent then $s_1s_2a=0\bigl(=s_1(s_2a)\bigr)$ implies $0=(s_2a)\tilde{s}_1=s_2(a\tilde{s}_1)$ with some
$\tilde{s}_1\in S$ and therefore $(a\tilde{s}_1)\tilde{s}_2=0$ with some $\tilde{s}_2\in S$. By multiplicativity,
$\tilde{s}_1\tilde{s}_2$ belongs to $S$.

 \vskip .2cm
Let $s=h_{\mu}+k$ and let $a$ be an element of $A$ such that $(h_{\mu}+k)a=0$. Let $a=a_1+...+a_N$ be the
$\mathfrak{h}$-weight decomposition of $a$ as in the check of the condition (\ref{loou1}). Then by the weight
argument, $(h_{\mu}+k)a_j=0$ for all $j=1,...,N$. But $(h_{\mu}+k)a_j=a_j(h_{\mu}+k+\mu_j)$ by homogeneity of $a_j$.
Therefore, for $\tilde{s}:=(h_{\mu}+k+\mu_1)(h_{\mu}+k+\mu_2)...(h_{\mu}+k+\mu_N)$ we obtain
$$a\tilde{s}=(a_1+...+a_N)(h_{\mu}+k+\mu_1)...(h_{\mu}+k+\mu_N)=0\ .$$
The check is completed. \hfill$\Box$

\section*{Acknowledgments}
 S. K. was 
supported by the RFBR grant 11-01-00980,  joint SFBRU-RFBR grant 11-02-90453, 
and by Federal Agency for Science and Innovations of Russian Federation under contract
14.740.11.0347.

\end{document}